\documentclass[12pt]{amsart}
\usepackage{amssymb}
\usepackage{amsmath}
\usepackage{amsthm}
\usepackage{mathrsfs}
\usepackage{amsfonts}

\usepackage[english]{babel}

\usepackage{parskip}

\usepackage{tikz,pgfplots}

\usepackage{rotating}

\usepackage{tikz-cd}

\usetikzlibrary{cd}

\usetikzlibrary{positioning}
\usetikzlibrary{arrows}
\usepackage{xcolor}

\usepackage{manfnt}

\title{On cardinality bounds for $\theta^n$-Urysohn spaces}
\date{}
\author{Fortunata Aurora Basile}
\address{University of Messina}
\email{basilef@unime.it}

\author{Nathan Carlson}
\address{California Lutheran University}
\email{ncarlson@callutheran.edu}

\author{Jack Porter}
\address{University of Kansas}
\email{porter@ku.edu}

\begin{document}

\maketitle

\newtheorem {theorem}{Theorem}
\newtheorem {lemma}{Lemma}
\newtheorem {cor}{Corollary}
\newtheorem {defin}{Definition}
\newtheorem{remark}{Remark}
\newtheorem{claim}{Claim}
\newtheorem{es}{Example}
\newtheorem{quest}{Question}
\newtheorem{pro}{Property}
\newtheorem{prop}{Proposition}

\begin{abstract}
We introduce the class of $\theta^{n}$-Urysohn spaces and the $n$-$\theta$-closure operator. $\theta^n$-Urysohn spaces generalize the notion of a Urysohn space. We estabilish bounds on the cardinality of these spaces and cardinality bounds if the space is additionally homogeneous.
\end{abstract}



\section{Introduction}
We shall follow notations from \cite{ENG} and \cite{H}. For a space $X$, we denote by $\chi(X)$, $\psi(X)$, $\pi_\chi(X), c(X), t(X)$ the \emph{character}, \emph{pseudocharacter}, $\pi$\emph{-character}, \emph{celluarity}, \emph{tightness} of a space $X$, respectively \cite{ENG}.\\
\indent Recall that a space $X$ is \textit{Urysohn} if for every two distinct points $x,\;y\in X$ there are open sets $U$ and $V$ such that $x\in U$, $y\in V$ and $\overline{U}\cap\overline{V}=\emptyset$. Related to Urysohn spaces we have the notion of the $\theta$\textit{-closure} \cite{VEL}. The $\theta$\textit{-closure} of a set $A$ in a space $X$ is the set $cl_{\theta}(A)=\{x\in X:$ for every neighborhood $U\ni x , \overline{U}\cap A\neq\emptyset\}$. $A$ is said to be $\theta$-\textit{closed} if $A = cl_\theta(A)$. The $\theta$\textit{-tightness} of $X$ at $x\in X$ is $t_{\theta}(x,X)=\min \{\kappa :$ for every $A\subseteq X$ with $x\in cl_{\theta}(A)$ there exists $B\subseteq A$ such that $|B|\leq \kappa$ and $x\in cl_{\theta}(B)\}$. The $\theta$\textit{-tightness of} $X$ is $t_{\theta}(X)=sup\{t_{\theta}(x,X):\;x\in X\}$ \cite{CK}. We have that tightness and $\theta$-tightness are independent (see Example 11 and Example 12 in \cite{PCC}), but if $X$ is a regular space then $t(X)=t_{\theta}(X)$. We say that a subset $A$ of $X$ is $\theta$\textit{-dense} in $X$ if $cl_{\theta}(A)=X$. The $\theta$-\textit{density} of $X$ is $d_{\theta}(X)=min\{\kappa:\;A\subseteq X\;\hbox{, }A\;\hbox{is a dense subset of }X\hbox{and }|A|\leq \kappa\}$. If $X$ is a Hausdorff space, the \textit{closed pseudocharacter of a point} $x$ in $X$ is $\psi_{c}(x,X)=\min\{|{\mathcal U}| : {\mathcal U}$ is a family of open neighborhoods of $x$ and $\left\{x\right\}$ is the intersection of the closure of $\mathcal U\}$. The \textit{closed pseudocharacter of} $X$ is $\psi_{c}(X)=sup\{\psi_{c}(x,X):\;x\in X\}$ (see \cite{SH} where it is called $S\psi(X)$). If $X$ is a Urysohn space, the $\theta$-\textit{pseudocharacter of a point} $x$ in $X$ is $\psi_{\theta}(x,X)=\min\{|{\mathcal U}| : {\mathcal U}$ is a family of open neighborhoods of $x$ and $\left\{x\right\}$ is the intersection of the $\theta$-closure of the closure of $\mathcal U\}$. The $\theta$-\textit{pseudocharacter of} $X$ is $\psi_{\theta}(X)=sup\{\psi_{\theta}(x,X):\;x\in X\}$ \cite{BBC}.\\
A collection $\mathcal V$ of open subsets of $X$ is called \emph{Urysohn-cellular}, if $O_1,\; O_2$ in $\mathcal V$ and $O_1\neq O_2$ implies
$\overline{O_1}\cap\overline{O_2}=\emptyset$.
The \emph{Urysohn-cellularity} of a space X is $Uc(X) =\ sup\{|{\mathcal V}| : {\mathcal V}$ is Urysohn-cellular$\}$. Of course, $Uc(X) \leq c(X)$. The \textit{almost Lindel\"of degree} of a subset $Y$ of a space $X$ is $aL(Y,X)=\min \{\kappa :$ for every cover $\mathcal{V}$ of $Y$ consisting of open subsets of $X$, there exists $\mathcal{V'}\subseteq\mathcal{V}$ such that $|\mathcal{V'}|\leq \kappa$ and
$\bigcup\{\overline{V}:\;V\in\mathcal{V'}\}=Y$\}.  The function $aL(X,X)$ is called the \textit{almost Lindel\"of degree} of $X$ and denoted by $aL(X)$ (see \cite{WD} and \cite{H}). The \textit{almost  Lindel\"of degree of $X$ with respect to closed subsets of $X$} is $aL_{c}(X)=\sup\{aL(C,X):\;C\subseteq X\;is\;closed\}$.
For a subset A of a space X we will denote by $[A]^{\leq\lambda}$ the
family of all subsets of A of cardinality $\leq\lambda$.

In Section \ref{S1} we give the definition of $n$-$\theta$-closure and using Example \ref{EN} we distinguish this operator from the $\theta^{n}$-closure defined in \cite{DG}. From the definition of $n$-$\theta$-closure, we are able to generalize the concept of Urysohn spaces introducing the $\theta^{n}$-Urysohn spaces and also we prove two characterizations (Proposition \ref{Car} and Proposition \ref{Car1}). We have that a $\theta^{(n+1)}$-Urysohn space is a $\theta^{n}$-Urysohn space and the Example \ref{EX3} shows that the converse is not true. With Example \ref{EsSn} we distinguish $\theta^{n}$-Urysohn spaces from $S(n)$-spaces defined in \cite{DG}. We note further that a $\theta^1$-Urysohn space is Urysohn.

In Section \ref{S2} we introduce new cardinal functions: the $n$-$\theta$-\textit{almost Lindel\"of degree} of a space $X$ (denoted by $\theta^{n}$-$aL(X)$), the $n$-$\gamma$-\textit{tightness} of a space $X$ (denoted by $t_{\gamma}^{n}(X)$), the $n$-$\gamma$-\textit{pseudocharacter} of the space $X$ (denoted by $\psi_{\gamma}^{n}(X)$), and the $\theta^{n}$-\textit{Urysohn cellularity} of the space $X$ (denoted by $\theta^{n}$-$Uc(X)$) in order to extend some known cardinality bounds for Urysohn spaces in the case of $\theta^{n}$-Urysohn spaces. In particular we prove the following:
\begin{itemize}
\item if $X$ is a $\theta^{n}$-Urysohn space, then $|X|\leq 2^{\psi_{\gamma}^{n}(X)t_{\gamma}^{n}(X)\theta\hbox{-}aL(X)}$ (Theorem \ref{3}). For $n=1$ we have Theorem 2 in \cite{BC}.
\item if $X$ is a $\theta^{n}$-Urysohn space, then $|X|\leq 2^{\theta^{n}\hbox{-}Uc(X)\chi(X)}$ (Theorem \ref{cardbound}). For $n=1$ we have Theorem 9 in \cite{Sc}.

\end{itemize}

Many cardinality bounds for general spaces have corresponding ``companion'' cardinality bounds for homogeneous topological spaces. In Theorem~\ref{homogeneousbound} in Section \ref{S3} we give the homogeneous companion bound to Theorem~\ref{cardbound}. In fact in Theorem \ref{PHbound} we prove that if $X$ is a power homogeneous and $\theta^{n}$-Urysohn space then $|X|\leq 2^{\theta^{n}\hbox{-}Uc(X)\pi_{\chi}(X)}$. This generalizes Theorem 13 in \cite{BCCS} and it is a modification of Theorem 2.3 in \cite{CR}.

\section{$n$-$\theta$-closure, $n$-$\gamma$-closure and $\theta^n$-Urysohn spaces}\label{S1}
In this section we introduce two closure operators and new axioms of separation.
\begin{defin}\rm
Let $X$ be a space. For $n\in\omega$, the \textit{n}-$\theta$-\textit{closure} of a subset $A$ of $X$ is
$$cl_{\theta}^{n}(A)=\underbrace{cl_{\theta}cl_{\theta}...cl_{\theta}}_{n\hbox{-}times}(A).$$

A subset $A$ of $X$ is called $n$-$\theta$-\textit{closed} if $A=cl_{\theta}^{n}(A)$.
\end{defin}

It is natural to compare this new closure operator with the following closure operator introduced in \cite{DG} by Dikranjan and Giuli.\\
\begin{defin}\rm\cite{DG}
Let $X$ be a topological space, $M$ a subset of $X$ and $n\in\omega$, $n\geq1$. We say that $x\in cl_{\theta^{n}}(M)$ if for every chain $U_{1}\subset U_{2}\subset...\subset U_{n}$ of open neighborhood of $x$ such that $\overline{U_{i}}\subset U_{i+1}$ for every $i=1,...,n\hbox{-}1$, $\overline{U_{n}}\cap M\neq\emptyset$. If $n=1$ we have the $\theta$-closure of a set.
\end{defin}

We give a relationship between $n$-$\theta$-closure and $\theta^{n}$-closure. The next proposition follows directly from \cite{DG}.
\begin{prop}\label{RC}
If $X$ is a space and $A$ is a subset of $X$, then $cl_{\theta}^{n}(A)\subseteq cl_{\theta^{n}}(A)$.
\end{prop}
If $U$ is an open subset of a space $X$ we have the following.
\begin{prop}\cite{S}
Let $X$ be a space and $U$ an open subset of $X$, then $cl_{\theta^{2}}(U)=cl_{\theta}(\overline{U})=cl_{\theta}^{2}(U)$.
\end{prop}

With Example \ref{EN} we can show that the containment in Proposition \ref{RC} can be strict even if $A$ is an open subset of $X$. First we need to recall the following example called the \textit{Tychonoff spiral} which will be used to prove some results.

\begin{es}\rm\label{TY}
The \textit{Tychonoff plank}, denoted by $T$,  is  the subspace
$(\omega_1+1)\times (\omega+1)\backslash\{(\omega_1,\omega)\}$ of the product space $(\omega_1+1)\times (\omega+1)$ where $\omega_1+1$ and $\omega +1$ are the usual ordered ordinals.

For $n \in \omega$, let  $T_n = T \times {n}$, $Y = \uplus_{i\in \omega} T_{i}$. For $n$ odd, identify the points $(\omega_1,k,n)$ with $(\omega_1,k,n+1)$ for $k \in \omega$ and 
for $n$ even, identify the points $(\alpha,\omega, n)$ with $(\alpha,\omega, n+1)$ when $\alpha \in \omega_1$.
Denote this new space $Z$.  $Z$ is called the \textit{Tychonoff spiral}.
\end{es}

\begin{es}\label{EN}
A space $X$ such that if $U$ is an open subset of $X$ then $cl_{\theta}^{n}(U)\neq cl_{\theta^{n}}(U)$.
\end{es}

For $k, n \in \omega$  and $ k \geq 3$, here is an example of a space $Y_k$ and open set $V \in \tau(Y_k)$ such that  $cl_{\theta}^n(V) = cl(V)$  and  $cl_{{\theta^n}}(V) =$

$\left \{\begin{array}{l@{\quad \quad}l} cl(V)&\text{if $n \leq k-1$ }\vspace{0mm}\\ 
\vspace{-1mm} cl(V) \cup \{p\}&\text{if $n >k-1$}\end{array}\right.$.
 
Let $Y_k$ be the subspace $T_1 \cup \cdots \cup  T_k$ of $Z$ in Example \ref{TY} plus an additional point $\{p\}$ with this topology: $U \in \tau(Y_k)$ if $U\backslash\{p\}$ is open in the subspace $T_1 \cup \cdots \cup T_k$, and $p \in U$ implies there are $\alpha \in \omega_1$ and $m \in \omega$ such that $(\alpha, \omega_1)  \times (m, \omega] \times \{1\} \subset U$.

\noindent The space $Y_3=T_1\cup T_2\cup T_3 \cup \{p\}$ looks like this:

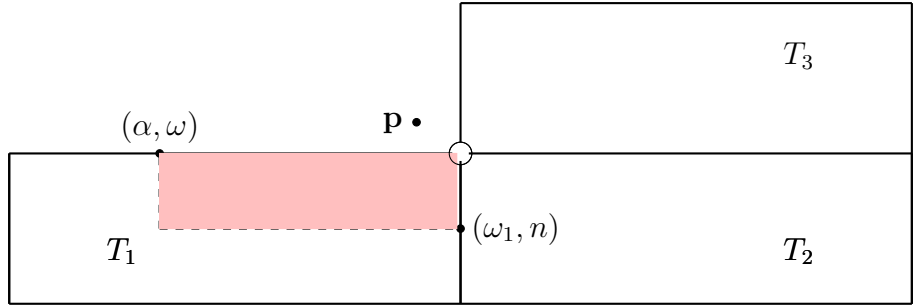
\begin{tikzpicture}

\draw[thick] (0,0) -- (12,0);
\draw[thick] (0,2) -- (5.89,2);
\draw[thick] (6.11,2) -- (12,2);
\draw[thick] (0,0) -- (0,2);
\draw[thick] (6,0) -- (6,1.9);
\draw[thick] (12,0) -- (12,2);

\draw (6,2) circle (1.5mm);
\node at (1.5, .7)   {$T_1$};
\node at (10.5, .7)   {$T_2$};

\node[above] at (2, 2)   {($\alpha, \omega$)};
\draw[fill] (2,2) circle [radius=0.05];
\draw[fill] (6,1) circle [radius=0.05];
\node[right] at (6, 1)   {($\omega_1, n$)};
\node at (5.1, 2.4)   {$\bf{p}$};
\draw[thick, dashed] (2,1) -- (6,1);
\draw[thick, dashed] (2,1) -- (2,2);

\draw[fill] (5.8,2.8, 1)  circle [radius=0.05];

\draw [fill=pink,pink] (2,1) rectangle (5.95,2); 
\draw[thick] (6,4) -- (12,4);
\draw[thick] (0,0) -- (0,2);
\draw[thick] (6,0) -- (6,1.85);
\draw[thick] (6,2.15) -- (6,4);
\draw[thick] (12,0) -- (12,4);
\draw (6,2) circle (1.5mm);
\node at (1.5, .7)   {$T_1$};
\node at (10.5, .7)   {$T_2$};
\node at (10.5, 3.3)   {$T_3$};
\end{tikzpicture}

Consider the open set $V = T_k\backslash (\{\omega_1\}\times\omega)$. We can easily seen that $cl^n_{{\theta}}(V) = cl(V)$ for $n \geq 1$; however,   $cl_{{\theta^n}}(V) = cl(V) \cup \{p\}$ for $n > k-1$ and $cl_{{\theta^n}}(V) = cl(V) $ for $n \leq k-1$.

Using the $n$-$\theta$-closure we introduce new axioms of separation as follows.

\begin{defin}\rm
We say that a space $X$ is $\theta^n$-$Urysohn$, for every $n\in\omega$, if for every $x,\;y\in X$ with $x\neq y$, there exist open subsets $U$ and $V$ of $X$ with $x\in U$ and $y\in V$ such that $cl_{\theta}^{n}(U)\cap cl_{\theta}^{n}(V)=\emptyset$.
\end{defin}

Observe that $\theta^1$-Urysohn spaces are Urysohn spaces. It is clear that a $\theta^{(n+1)}$-Urysohn space is a $\theta^n$-Urysohn space for every $n\in\omega$, but the converse is not true as the following example shows.

\begin{es}\label{EX3}
A $\theta^m$-Urysohn not $\theta^{(m+1)}$-Urysohn space.
\end{es} 

We present a  $\theta^m$-Urysohn space $X_{4m+1}$ that is not $\theta^{(m+1)}$-Urysohn where $m \in \omega$. Let $\mathbb {R}= \bigcup_{n \in \omega} A_n$ where ${A_n}'s$ are  pairwise disjoint, each $A_n$ is dense in $\mathbb{R} $, and $0 \in A_1$.  Let  $A_{4m+1}' = A_{4m+1} \cup \{0'\}$ where $0' \not\in \mathbb{R}$.  Let $X_{4m+1} = \bigcup_{n=1}^{4m}A_n \cup \{A'_{4m+1}\}$. 

\noindent For $a,b \in \mathbb{R}$ where $a < b$, let 

$(a,b)' =  \left \{\begin{array}{l@{\quad \quad}l} (a,b)&\text{if $0 \notin (a,b)$ }\vspace{0mm}\\ 
\vspace{-1mm} ((a,b)\backslash\{0\}) \cup\{0'\}&\text{if $0 \in (a,b)$}\end{array}\right.$  
 \vskip 2mm

If $a, b \in \mathbb{R}$ and $a < b$, an open base for $X_{4m+1}$ is generated by these families of sets:  

\noindent(1)  if $n\in\omega$ is odd and $1 \leq n \leq 4m-1$, $(a,b) \cap A_n$ is open, \newline 
(2)  if $n\in\omega$ is even, $2 \leq n \leq 4m$, $a,b\in\mathbb{R}$ and $0 \notin (a,b)$,  $(a,b)\cap (A_{n-1} \cup A_n \cup A_{n+1})$ is open, and \newline
(3) if $n = 4m+1$, $(a,b)'\cap A'_{4m+1}$ is open.

The following picture represents the space $X_9$.

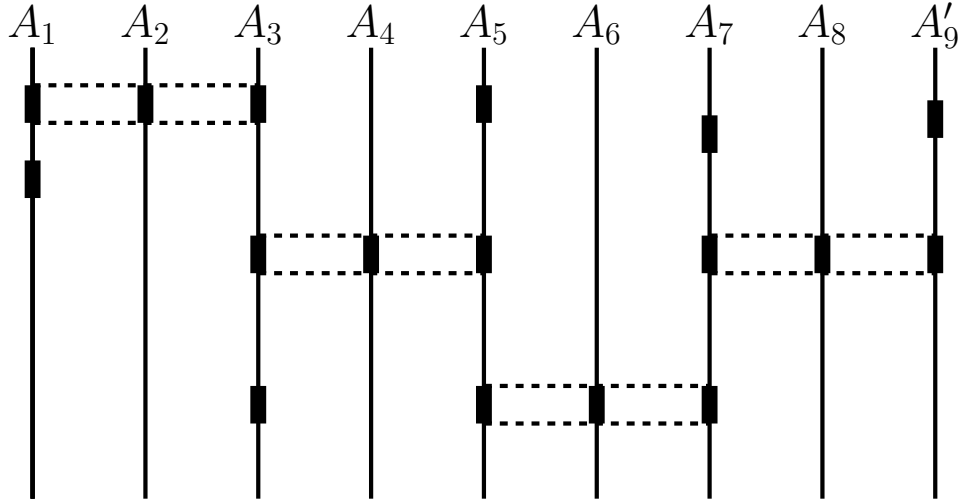
\begin{tikzpicture}
\draw[ultra thick] (0,0) --(0,6);
\draw[ultra thick] (1.5,0) --(1.5,6);
\draw[ultra thick] (3,0) --(3,6);
\draw[ultra thick] (4.5,0) --(4.5,6);
\draw[ultra thick] (6,0) --(6,6);
\draw[ultra thick] (7.5,0) --(7.5,6);
\draw[ultra thick] (9,0) --(9,6);
\draw[ultra thick] (10.5,0) --(10.5,6);
\draw[ultra thick] (12,0) --(12,6);

\draw [dashed, ultra thick] (0,5) -- (3,5);
\draw [dashed, ultra thick] (0,5.5) -- (3,5.5);

\draw [line width=6](0,5)--(0,5.5);
\draw [line width=6](1.5,5)--(1.5,5.5);
\draw [line width=6](3,5)--(3,5.5);

\draw [line width=6](7.5,1)--(7.5,1.5);
\draw [line width=6](9,3)--(9,3.5);
\draw [line width=6](10.5,3)--(10.5,3.5);
\draw [line width=6](12,3)--(12,3.5);

\draw [line width=6](4.5,3)--(4.5,3.5);

\draw [line width=6](6,1)--(6,1.5);
\draw [line width=6](6,5)--(6,5.5);
\draw [line width=6](9,4.6)--(9,5.1);
\draw [line width=6](3,1)--(3,1.5);

\draw [dashed, ultra thick] (3,3.5) -- (6,3.5);
\draw [dashed, ultra thick] (3,3) -- (6,3);

\draw [dashed, ultra thick] (9,3.5) -- (12,3.5);
\draw [dashed, ultra thick] (9,3) -- (12,3);


\draw [dashed, ultra thick] (6,1.5) -- (9,1.5);
\draw [dashed, ultra thick] (6,1) -- (9,1);

\draw [line width=6](0,4)--(0,4.5);
\draw [line width=6](3,3)--(3,3.5);

\draw [line width=6](6,3)--(6,3.5);
\draw [line width=6](12,4.8)--(12,5.3);

\draw [line width=6](9,1)--(9,1.5);

\node at (0,6.3) {\Large{$A_1$}};
\node at (1.5,6.3) {\Large $A_2$};
\node at (3,6.3) {\Large $A_3$};
\node at (4.5,6.3) {\Large $A_4$};
\node at (6,6.3) {\Large $A_5$};
\node at (7.5,6.3) {\Large $A_6$};
\node at (9,6.3) {\Large $A_7$};
\node at (10.5,6.3) {\Large $A_8$};
\draw[ultra thick] (0,0) --(0,6);
\draw[ultra thick] (0,0) --(0,6);
\node at (12,6.3) {\Large $A'_9$};
\end{tikzpicture}

\vskip 2mm
\noindent Let $a,b,c,d\in\mathbb{R}$, $U = (a,b)\cap A_1$. Then $cl_{\theta}(U) = cl(U) \subseteq [a,b]\cap (A_1 \cup A_2)$.  Let  $V = (c,d)\cap(A_3 \cup A_4 \cup A_5)$.  Then $cl_{\theta}(V) = cl(V) = [c,d] \cap (A_2 \cup A_3 \cup A_4 \cup A_5 \cup A_6)$.  It follows that $cl^2_{\theta}(U) = [a,b]\cap (A_1 \cup A_2 \cup A_3 \cup A_4)$.  By induction, $cl^m_{\theta}(U) = [a,b]\cap (A_1 \cup A_2 \cup \cdots   \cup A_{2m})$.  Likewise, starting from the right-hand subspace $A'_{4m+1}$ with $U = (a,b)'\cap A_{4m+1}$, we have $cl^m_{\theta}(U) = [a,b]'\cap (A'_{4m+1} \cup A_{4m} \cup \cdots   \cup A_{2m+1})$.\\
We have the following consequences:

(I): Every pair of points in $X_{4m+1}$, except $0$ and $0'$, are $\theta^n$-Urysohn separated for all $n \in \omega$. 

(II):  Let $m \geq 1$ and $a,b,c,d,\in\mathbb{R}$. Let $0 \in (a,b),\; 0' \in (c,d)',\; U = (a,b) \cap (c,d).$  Then  $cl^m_{\theta}(U)  = [a,b]\cap (A_1 \cup A_2 \cup \cdots   \cup A_{2m})$ and $cl^m_{\theta}(U') = [a,b]'\cap (A_{4m+1} \cup A_{4m} \cup \cdots   \cup A_{2m+2})$.  Thus, $cl^m_{\theta}(U) \cap cl^m_{\theta}(U') = \varnothing$ and $X_{4m+1}$ is $\theta^m$-Urysohn.  On the other hand,  $cl^{m+1}_{\theta}(U) \cap cl^{m+1}_{\theta}(U') \supseteq [a,b] \cap [c,d] \cap A_{2m+1} \ne \varnothing$; so, $X_{4m+1}$ is not $\theta^{(m+1)}$-Urysohn.  

(III):  Using the argument in (II), $X_{5}$ is a Urysohn not firmly Urysohn space (a space $X$ is called \textit{firmly Urysohn} \cite{BaC} if for every $x,y\in X$ with $x\neq y$ there exist open subsets $U,V$ of $X$ with $x\in U$ and $y\in V$ such that $\overline{U}\cap cl_{\theta}(\overline{V})=\emptyset$). In fact, for $0 \in (a,b), 0' \in (c,d)'$ and $U = (a,b) \cap (c,d)$, we have $cl^2_{\theta}(U) = [a,b]\cap (A_1 \cup A_2 \cup A_3 \cup A_4)$ and  $\overline{U'}=cl^1_{\theta}(U') = [a,b]'\cap (A_{4} \cup A_{5})$.  Thus, $cl^2_{\theta}(U) \cap 
cl^1_{\theta}(U') \supseteq [a,b] \cap [c,d] \cap A_{4} \neq\emptyset$. This answers Question 2.1 in \cite{BaC}.

With Example \ref{EN} we distinguish the $\theta^{n}$-closure from the $n$-$\theta$-closure. It is natural to investigate the relation between $S(n)$-spaces defined in \cite{DG} and $\theta^{n}$-Urysohn spaces.

\begin{defin}\cite{DG}\rm
A space $X$ is called $S(n)$-space if for every $x,\;y\in X$, $x\neq y$, are $\theta^{n}$-separated. A point $x\in X$ is $\theta^{n}$-\textit{separated from a subset} $M$ of $X$ if $x\notin cl_{\theta^{n}}(M)$. For $n>0$, the relation \textit{being} $\theta^{n}$-\textit{separated} between points is symmetric.
\end{defin}

\begin{es}\label{EsSn}
A $\theta^n$-Urysohn not $S(n)$-space.
\end{es}
\proof
In the Tychonoff spiral $Z$ described in Example \ref{TY}, let $X_n$ denote the subspace $T_1 \cup T_2 \cup \cdots \cup T_n$ plus two additional points $\{p,q\}$ with this topology: $U \in \tau(X_n)$ if $U\backslash\{p, q\}$ is open in the subspace $T_1 \cup T_2 \cup \cdots T_n$, $p \in U$ implies there are $\alpha \in \omega_1$ and $n \in \omega$ such that $(\alpha, \omega_1) \times (n, \omega] \times \{1\} \subset U$, and $q \in U$  implies there are $\alpha \in \omega_1$ and $n \in \omega$ such that $(\alpha, \omega_1) \times (n, \omega] \times \{n\} \subset U$.

The spaces $X_3$ and $X_4$ are $\theta^n$-Urysohn for $n \in \omega$ and $S(1)$ but not $S(2)$.
For $k \in \omega$,  the spaces $X_{2k+1}$ and $X_{2k+2}$ are $\theta^n$-Urysohn for $n \in \omega$ and $S(k)$ but not $S(k+1)$.
\endproof

\begin{quest}\rm
Does there exist a $S(n)$ not $\theta^{n}$-Urysohn space for every $n\geq2 \in \omega$?
\end{quest}

In order to prove a characterization for $\theta^{n}$-Urysohn spaces, we need to introduce a new operator called $\gamma$-$n$-closure.

\begin{defin}\rm
Let $X$ be a space, the \textit{n}-$\gamma$-\textit{closure} of a subset $A$ of $X$, for every $n\in\omega$ is
$$cl_{\gamma}^{n}(A)=\{x\in X:\;\hbox{for every open subset}\;U\;\hbox{of}\;X,\;cl_{\theta}^{n}(U)\cap A\neq\emptyset\}$$

\noindent A subset $A$ of $X$ is called $n$-$\gamma$-closed if $A=cl_{\gamma}^{n}(A)$.
\end{defin}

If $n=1$ we have the $\theta$-closure and if $n=2$ we have the $\gamma$-closure defined in \cite{BBC}. We have that $cl_{\gamma}^{n}(A)\supseteq cl_{\theta}^{n}(A)$ for every subset $A$ of $X$ and with the next example we can show that this inequality can be strict.
\begin{es}\rm
If $X$ is the Bing's Tripod space then for every $x\in X$ $cl_{\theta}^{n}(\{x\})=\{x\}$ and $cl_{\gamma}^{n}(\{x\})=X$.
\end{es}

\begin{prop}\label{Car}
A space $X$ is $\theta^n$-Urysohn if and only if for every $x\in X$ and for every family ${\mathcal U}_{x}$ of open neighborhoods of $x$, $\{x\}=\bigcap_{U\in{\mathcal U}_{x}}cl_{\gamma}^{n}(cl_{\theta}^{n}(U))$.
\end{prop}
\proof
Let $X$ be $\theta^n$-Urysohn and $x\in X$. $\forall y\in X\setminus\{x\}$, there exists $U_{y},\;V_{y}$ open subsets of $X$ with $x\in U_{y}$ and $y\in V_{y}$ such that $cl_{\theta}^{n}(U)\cap cl_{\theta}^{n}(V)=\emptyset$. This means $y\notin cl_{\gamma}^{n}(cl_{\theta}^{n}(U_{y}))$ and $\{x\}=\bigcap_{y\in X\setminus\{x\}}cl_{\gamma}^{n}(cl_{\theta}^{n}(U_{y}))$.

For the converse, let $x,\;y\in X$ with $x\neq y$. Then there exists an open subset $V$ of $X$ containing $x$ such that $y\notin cl_{\gamma}^{n}(cl_{\theta}^{n}(U))$ so there exists an open subset $U$ of $X$ containing $y$ such that $cl_{\theta}^{n}(U)\cap cl_{\theta}^{n}(V)=\emptyset$. This means $X$ is $\theta^n$-Urysohn.
\endproof

We want to show that $n$-$\theta$-closure is productive. The next proposition shows that 1-$\theta$-closure is productive and the proof can be easily extended to a the $n$-$\theta$-closure operator for every $n\in\omega$.

\begin{prop}\label{P30}
If $\{X_{s}\}_{s\in S}$ is a family of spaces, $X=\prod_{s\in S}X_{s}$ and $A_{s}$ is a subset of $X_{s}$ for every $s\in S$ then $cl_{\theta}(\prod_{s\in S}A_{s})=\prod_{s\in S}cl_{\theta}(A_{s})$.
\end{prop}
\proof
Let $x\in cl_{\theta}(\prod_{s\in S}A_{s})$ iff for every member $\prod_{s\in S} W_{s}$ of the canonical base containing $x$ then $\emptyset\neq \overline{\prod_{s\in S} W_{s}}\cap \prod_{s\in S}A_{s}= \prod_{s\in S}\overline{W_{s}}\cap\prod_{s\in S}A_{s}=\prod_{s\in S}\overline{W_{s}}\cap A_{s}$ iff $x\in \prod_{s\in S}cl_{\theta}(A_{s})$.
\endproof

The following is another characterization of $\theta^{n}$-Urysohn spaces. We know that a space $X$ is Hausdorff iff the diagonal $\Delta_{X}$ is closed in $X\times X$ and that a space $X$ is Urysohn iff the diagonal $\Delta_{X}$ is $\theta$-closed in $X\times X$. We have the same result for $\theta^{n}$-Urysohn spaces.

\begin{prop}\label{Car1}
A space $X$ is $\theta^{n}$-Urysohn iff the diagonal $\Delta_{X}$ is $n$-$\gamma$-closed in $X\times X$.
\end{prop}
\proof
$\Delta_{X}$ is $n$-$\gamma$-closed in $X\times X$. Suppose $x,\;y\in X$ are distinct points then $(x,y)\notin \Delta_{X}$ iff there exists a basic open set $B$ in $X\times X$ containing $(x,y)$ such that $cl_{\theta}^{n}(B)\cap \Delta_{X}=\emptyset$. The basic open set $B$ is expressed as the product of $U\times V$ where $U$ is an open set in $X$ containing $x$ and $V$ is an open set in $X$ containing $y$. But $cl_{\theta}^{n}(B)=cl_{\theta}^{n}(U)\times cl_{\theta}^{n}(V)$ and this means $x\notin cl_{\gamma}^{n}(cl_{\theta}^{n}V)$ iff $X$ is $\theta^{n}$-Urysohn (by Proposition \ref{Car}).
\endproof

\section{Cardinality bounds involving $\theta^n$-Urysohn spaces}\label{S2}

In 1988, Bella and Cammaroto \cite{BC} proved that if $X$ is a Urysohn space, then $|X|\leq 2^{\chi(X) aL(X)}$. Hodel mentioned in his survey \cite{H} this variation: If $X$ is a Urysohn space then $|X|\leq 2^{t(X)\psi(X)aL_{c}(X)}$.

First some definitions.

\begin{defin}\rm
If $X$ is a space and $x\in X$, the $n\hbox{-}\gamma\hbox{-}tightness$ of $x$ with respect to $X$ is $t_{\gamma}^{n}(x,X)=min\{k:\;\forall x\in cl^{n}_{\gamma}(A),\;\exists B\in[A]^{\leq k}:\;x\in cl^{n}_{\gamma}(B)\}$. The $n\hbox{-}\gamma\hbox{-}tightness$ of $X$ is $t_{\gamma}^{n}(X)=sup_{x\in X}t_{\gamma}^{n}(x,X)$.

\end{defin}

If $n=1$, the $1$-$\gamma$-tightness of a space $X$ is the $\theta$-tightness of a space $X$ defined in \cite{CK}.

\begin{prop}\label{R1}
If $X$ is a space such that $t_{\gamma}^{n}(X)\leq \kappa$ and if $H\subseteq X$ such that $H=\cup\{cl_{\gamma}^{n}(A):\;A\in[H]^{\leq\kappa}$\}, then  $H$ is $n$-$\gamma$-closed.
\end{prop}
\proof
Let $t_{\gamma}^{n}(X)\leq \kappa$ and $H$ as in the statement. We want to prove that $cl_{\gamma}^{n}(H)\subseteq H$. Let $x\in cl_{\gamma}^{n}(H)$. Then there is $A\in [H]^{\leq\kappa}$ and $x\in cl_{\gamma}^{n}(A)$ but $cl_{\gamma}^{n}(A)\subseteq H$ and this means $x\in H$.
\endproof

If $n=1$ we have Proposition 1.2 in \cite{CK}. Using Proposition \ref{Car} we are able to generalize the notion of pseudocharacter for $\theta^{n}$-Urysohn spaces.

\begin{defin}\rm
If $X$ is a $\theta^n$-Urysohn space, the $n$-$\gamma$\textit{-pseudocharacter of a point} $x \in X$, denote by  $\psi^{n}_{\gamma}(x,X)$ is:

$$
\psi^{n}_{\gamma}(x,X)=min\{\kappa:\;\hbox{there is a family}\;{\mathcal B}\;\hbox{of open neighborhoods of}\;x\;:|{\mathcal B}|\leq\kappa$$ $$\hbox{and}\;\{x\}=\bigcap_{U\in{\mathcal B}}cl_{\gamma}^{n}(cl_{\theta}^{n}(U))\}. 
$$

The $n$-$\gamma$\textit{-pseudocharacter of} $X$ is $\psi^{n}_{\gamma}(X)=sup\{\psi^{n}_{\gamma}(x,X):\;x\in X\}$.
\end{defin}
If $n=1$, the $1$-$\gamma$-pseudocharacter of a space $X$ is $\theta$-pseudocharacter of a space $X$  defined in \cite{BBC}.

\noindent We can easily see that $\psi_{\gamma}^{n}(X)\leq \chi(X)$ for every $n\in\omega$.

The following represents the relation between the $n$-$\gamma$-tightness and the character.

\begin{prop}\label{4}
For every space $X$, $t_{\gamma}^{n}(X)\leq \chi(X)$.
\end{prop}
\proof
Let $x\in X$, $A$ a subset of $X$ such that $x\in cl_{\gamma}^{n}(A)$, and ${\mathcal V}_{x}$ an open neighborhood system of $x$ with $|{\mathcal V}_{x}|\leq\chi(x,X)$. $x\in cl_{\gamma}^{n}(A)$ so for every $V\in {\mathcal V}_{x}$, $cl_{\theta}^{n}(V)\cap A\neq \emptyset$. We choose $y_{V}\in cl_{\theta}^{n}(V)\cap A$ for every $V\in{\mathcal V}_{x}$ and put $B=\{y_{V}:\;V\in{\mathcal V}_{x}\}$. $B\subseteq A$, $x\in cl_{\gamma}^{n}(B)$ and $|B|\leq\chi(x,X)$. This proves $t_{\gamma}^{n}(x,X)\leq \chi(x,X)$.
\endproof

The definition below is the generalization of the almost Lindel\"of degree of a space.

\begin{defin}\rm
The $n\hbox{-}\theta$\textit{-almost Lindel\"of degree} of a subset $Y$ of a space $X$, for every $n\geq2\in\omega$ is

$\theta^{n}\hbox{-}aL(Y,X)=\min \{\kappa : $ for every cover $\mathcal{V}$ of $Y$ consisting of open subsets of  $X$,
there exists $\mathcal{V'}\subseteq\mathcal{V}$  such that $|\mathcal{V'}|\leq \kappa$ and $\bigcup\{cl_{\theta}^{n}(V):\;V\in\mathcal{V'}\}=Y\}.$

The function $\theta^{n}$-$aL(X,X)$, $n\geq2\in\omega$, is called $n\hbox{-}\theta$\textit{-almost Lindel\"of degree of the space} $X$ and denoted by $\theta^{n}$-$aL(X)$.
\end{defin}

$\theta^{1}$-$aL(X)$ is $aL(X)$ defined in \cite{WD} and $\theta^{2}$-$aL(X)$ is $\theta$-$aL(X)$ defined in \cite{BBC}.

\begin{prop}
For every space $X$, $\theta^{n+1}\hbox{-}aL(X)\leq \theta^{n}\hbox{-}aL(X).$
\end{prop}

In general we have that the almost Lindel\"of degree is not a hereditary cardinal function but it is hereditary with respect to $\theta$-closed subsets. We prove that the $\theta^{n}$-almost Lindel\"of degree is hereditary with respect to $n$-$\gamma$-closed sets.

\begin{prop}\label{PR}
If $C$ is a $n$-$\gamma$-closed subset of $X$ then $\theta^{n}$-$aL(C,X)\leq \theta^{n}$-$aL(X)$.
\end{prop}
\proof
Let $X$ be a topological space such that $\theta^{n}\hbox{-}aL(X)\leq \kappa$ and let $C\subseteq X$ be $n$-$\gamma$-closed set. $\forall x\in X\setminus C$ we have that there exists an open neighborhood $U_{x}$ of $x$ such that $cl_{\theta}^{n}(U_x)\subseteq X\setminus C$. Let $\mathcal{U}$ be a cover of $C$ consisting of open subsets of $X$. Then $\mathcal{V}=\mathcal{U}\bigcup\{U_{x}:\;x\in X\setminus C\}$ is an open cover of $X$ and since $\theta^{n}\hbox{-}aL(X)\leq \kappa$, there exists
$\mathcal{V'}\in[\mathcal{V}]^{\leq \kappa}$ such that $X=\bigcup \{cl_{\theta}^{n}(V): V\in {\mathcal{V'}}\}$. Then there exists
$\mathcal{V''}\in[\mathcal{U}]^{\leq \kappa}$ such that $C\subseteq \bigcup \{cl_{\theta}^{n}(V): V\in {\mathcal{V''}}\}$; this proves that
$\theta^{n}\hbox{-}aL(C,X)\leq \kappa$.
\endproof

We now give a bound for the cardinality of the $n$-$\gamma$-closure of a subset of a space. This bound allows us to obtain cardinality bounds for $\theta^{n}$-Urysohn spaces.

\begin{prop}\label{PPR}
If $X$ is $\theta^n$-Urysohn and $A \subseteq X$,  then $|cl_{\gamma}^n(A)| \leq 2^{\psi_{\gamma}^n(X) t_{\gamma}^n(X)}$.
\end{prop}
\proof
Let $\kappa = \psi_{\gamma}^n(X) t_{\gamma}^n(X)$.  For each $x \in X$, there is a family $\mathcal{B}_x $ of open sets in $X$ containing $x$ such that $\bigcap_{U\in{\mathcal{B}_x}} cl_{\gamma}^n(cl_{\theta}^n(U)) = \{x\}$ and $|\mathcal{B}_x| \leq \kappa$. Let $x \in cl_{\gamma}^n(A)$ and $V$ be an open set containing $\{x\}$.  As  $x \in cl_{\gamma}^n(A)$, $cl_{\theta}^n(U \cap V) \cap A \ne \varnothing$.  Thus, $cl_{\theta}^n(U) \cap cl_{\theta}^n(V) \cap A \ne \varnothing$ and it follows that $x \in cl_{\gamma}^n(cl_{\theta}^n(U) \cap A) $ for all $U \in \mathcal{B}_x$. As $t_{\gamma}^n(X) \leq \kappa$, for each $U \in \mathcal{B}_x$, there is $A_U \subseteq cl_{\theta}^{n}(U)\cap A$ such that $|A_U| \leq \kappa$ and $x \in  cl_{\gamma}^n(A_U)$.  Thus, $\{x\} \subseteq \bigcap_{U\in{\mathcal{B}_x}} cl_{\gamma}^n(A_U) \subseteq\bigcap_{U\in{\mathcal{B}_x}} cl_{\gamma}^n(cl_{\theta}^n(U \cap A)) \subseteq \bigcap_{U\in{\mathcal{B}_x}} cl_{\gamma}^n(cl_{\theta}^n(U)) = \{x\}$.  Now, $\{A_U:U \in \mathcal{B}_x\} \in [[A]^{\leq \kappa}]^{\leq \kappa}$ and $|cl_{\gamma}^n(A)| \leq |A|^{\kappa}$. 
\endproof

\begin{defin}\rm
We say that a subset $A$ of $X$ is $n$-$\gamma$\textit{-dense} in $X$ if $cl_{\gamma}^{n}(A)=X$.
The $n$-$\gamma$-density of $X$ is:
$$d_{\gamma}^{n}(X)=min\{\kappa:\;A\subseteq X\;\hbox{, }A\;\hbox{is a }\;n\hbox{-}\gamma\hbox{ dense subset of }X\hbox{and }|A|\leq \kappa\}.$$ 
\end{defin}

\begin{cor}
If $X$ is a $\theta^{n}$-Urysohn space then $|X|\leq d_{\gamma}^{n}(X)^{\psi_{\gamma}^n(X) t_{\gamma}^n(X)}$.
\end{cor}
\proof
Let $A$ be a $n$-$\gamma$-dense subset of $X$, i.e. $cl_{\gamma}^{n}(A)=X$, with $|A|= d_{\gamma}^{n}(X)$. From the above theorem we have that $|cl_{\gamma}^{n}(A)|\leq |A|^{\psi_{\gamma}^n(X) t_{\gamma}^n(X)}$, so $|X|\leq d_{\gamma}^{n}(X)^{\psi_{\gamma}^n(X) t_{\gamma}^n(X)}$.
\endproof

Next we will use a variation of a result by Hodel \cite{H} that is proved in \cite{CCP}.

\begin{theorem}\label{1}
Let $X$ be a set, $\kappa$ an infinite cardinal, $c:[X]^{\leq \kappa}\rightarrow \mathcal{P}(X)$ a function, and for $x \in X$, $\{V(x,\alpha) < \kappa\}$ a collection of subsets of $X$ with these three properties:

{\bf(ME)} for $A, B \in [X]^{\leq \kappa}, A \subseteq c(A),$ and $c(A) \subseteq c(B)$,

{\bf(C)} for $A \in [X]^{\leq \kappa}, |c(A)| \leq 2^{\kappa}$, and 

{\bf (C-S)} if $H \ne \varnothing, |H| \leq 2^{\kappa}, c(B)  \subseteq H$ for $B \in [H]^{\leq \kappa}$ and $q \notin H$,  then there is $A \in [H]^{\leq \kappa}$ and a function $f:A\rightarrow \kappa$ such that  $H \subseteq  \bigcup_{a \in A}V(x,f(x)) \subseteq X\backslash\{q\}. $

Then $|X| \leq 2 ^{\kappa}.$
\end{theorem}

\begin{theorem}\label{3}
Suppose $X$ is $\theta^n$-Urysohn.  Then $|X| \leq 2^{\psi_{\gamma}^n(X)t_{\gamma}^n(X)\theta\hbox{-}aL^n(X)}$.
\end{theorem}
\proof
Let $\kappa = \psi_{\gamma}^n(X)t_{\gamma}^n(X)\theta\hbox{-}aL^n(X)$.  We will apply Theorem \ref{1} by verifying the properties (ME), (C), and (C-S).  For $A \in [X]^{\leq \kappa}$, define $c(A) = cl_{\gamma}^n(A)$; it is easy to check that the function ``c" satisfies (ME). By Proposition \ref{PPR}, (C) is satisfied. 
In preparation to show (C-S), we first define ``$V(x,\alpha)$".  As $\psi_{\gamma}^n(X) \leq \kappa$, for each $x \in X$, there is a family $\mathcal{B}_x = \{U(x,\alpha): \alpha < \kappa\}$ of open sets in $X$ containing $x$ such that $\bigcap_{\alpha < \kappa} cl_{\gamma}^n(cl_{\theta}^n(U(x,\alpha))) = \{x\}$.  Let $V(x,\alpha) = cl_{\gamma}^{n}(cl_{\theta}^n(U(x,\alpha)))$. Suppose $\varnothing \ne H \subseteq X$ such that $c(B)  \subseteq H$ for $B \in [H]^{\leq \kappa}$ and $q \notin H$.  For each $p \in H$, there is $f(p) < \kappa$ such that $q \notin cl_{\gamma}^n(cl_{\theta}^n(U(x,f(p)))$.  Thus, $\{U(p,f(p)): p \in H\}$ is an open cover of $H$. $H$ is $n$-$\gamma$-closed by Proposition \ref{R1} so we can apply Proposition \ref{PR}. Therefore $\theta\hbox{-}aL^n(H,X) \leq \theta\hbox{-}aL^n(X) \leq \kappa$, and there is a subset $A \subseteq H$ such that $|A| \leq \kappa$ and $H \subseteq \bigcup_{p \in A} cl_{\theta}^n(U(x,f(p))) \subseteq \bigcup_{p \in A} cl_{\gamma}^n(cl_{\theta}^n(U(x,f(p)))) = \bigcup_{p \in A}V(p,f(p)) \subseteq X\backslash \{q\}$.
\endproof

\begin{cor}
If $X$ is a $\theta^n$-Urysohn space then $|X|\leq 2^{\theta^{n}\hbox{-}aL(X)\chi(X)}$.
\end{cor}
\proof
We have $\psi_{\gamma}^{n}(X)\leq\chi (X)$ and by Proposition \ref{4} we have $t_{\gamma}^{n}(X)\leq\chi(X)$. Applying Theorem \ref{3} we have $|X|\leq 2^{\theta\hbox{-}aL(X)\chi(X)}$.
\endproof
If $n=1$ we have the Bella-Cammaroto inequality \cite{BC}:
\begin{cor}
If $X$ is a Urysohn space then $|X|\leq 2^{aL(X)\chi(X)}$. 
\end{cor}

In \cite{Sc} Schr\"oder proved that if $X$ is a Urysohn space then $|X|\leq 2^{Uc(X)\chi(X)}$. We can generalize this result, following a similar proof, in the class of $\theta^n$-Urysohn spaces.

\begin{defin}\rm
Let $X$ be a space. A family of open subsets ${\mathcal V}$ of $X$ is called $\theta^n\hbox{-}Urysohn\;cellular$ if for every $V_{1},\;V_{2}\in{\mathcal V}$, $cl_{\theta}^{n}(V_{1})\cap cl_{\theta}^{n}(V_{2})=\emptyset$.
\end{defin}

\begin{defin}\rm
The $\theta^n\hbox{-}Urysohn\;cellularity$ of a space $X$ is 
$$\theta^{n}\hbox{-} Uc(X)=sup\{|{\mathcal V}|:\;{\mathcal V}\;\hbox{is}\;\theta^n\hbox{-}\hbox{Urysohn cellular}\}.$$
\end{defin}

If $n=1$, $\theta^1\hbox{-}Urysohn\;cellularity$ is the Urysohn cellularity defined in \cite{Sc}.

\begin{lemma}\label{L1}
Let $X$ be a space and $\theta^{n}\hbox{-}Uc(X)\leq\kappa$. Let $\{U_{\alpha}\}_{\alpha\in A}$ be a family of open subsets of $X$. Then there exists $B\subseteq A$ such that $|B|\leq\kappa$ and $\bigcup_{\alpha\in A} U_{\alpha}\subseteq cl_{\gamma}^{n}(\bigcup_{\beta\in B}cl_{\theta}^{n}(U_{\beta}))$.
\end{lemma}
\proof
Let $\mathcal V$ be the collection of all the open subsets of $X$ contained in some $U_{\alpha}$. By Zorn's Lemma, there exists a maximal $\theta^n$-Urysohn cellular family ${\mathcal W}\subset{\mathcal V}$. For every $W\in{\mathcal W}$ take $U_{\beta}\in\{U_{\alpha}:\;\alpha\in A\}$ such that $W\subseteq U_{\beta}$. We may assume $\beta\in B\subseteq A$ and $|B|\leq\kappa$. Assume $\bigcup_{\alpha\in A} U_{\alpha}\nsubseteq cl_{\gamma}^{n}(\bigcup_{\beta\in B}cl_{\theta}^{n}(U_{\beta}))$. Then there exists $\alpha_{0}\in A$, $x\in U_{\alpha_{0}}$ and $U_{x}$ open neighborhood of $x$ such that $cl_{\theta}^{n}(U_{x})\cap \bigcup_{\beta\in B}cl_{\theta}^{n}(U_{\beta})=\emptyset$. Then $U_{x}\cap U_{\alpha_{0}}\in{\mathcal V}$ and $cl_{\theta}^{n}(U_{x}\cap U_{\alpha_{0}})\cap \bigcup_{W\in{\mathcal W}}cl_{\theta}^{n}(W)=\emptyset$, contradicting the maximality of ${\mathcal W}$.
\endproof

We are now ready to generalize the result by Schr\"oder. 

\begin{theorem}\label{cardbound}
If $X$ is a $\theta^n$-Urysohn space then $|X|\leq 2^{\theta^{n}\hbox{-}Uc(X)\chi(X)}$.
\end{theorem}
\proof
Let $\kappa=\theta^{n}$-$Uc(X)\chi(X)$. For every $x\in X$ let ${\mathcal B}(x)$ an open neighborhood base of $x$ with $|{\mathcal B}(x)|\leq\kappa$. Construct an increasing sequence $\{C_{\alpha}\}_{\alpha<\kappa}$ of subsets of $X$ and a sequence $\{{\mathcal V}_{\alpha}\}$ of families of open subsets of $X$ such that:
\begin{enumerate}
\item $|C_{\alpha}|\leq\kappa$ for every $\alpha\leq\kappa$;
\item ${\mathcal V}_{\alpha}=\bigcup\{{\mathcal B}(c):\;c\in\bigcup_{\tau<\alpha} C_{\tau}\}$, $\alpha<\kappa^{+}$;
\item if $\{G_{\beta}:\beta<\kappa\}$ is a collection of subsets of $X$ and each $G_{\beta}$ is the union of the $n$-$\theta$-closures of $\leq\kappa$ many elements of ${\mathcal V}_{\alpha}$ and $\bigcup _{\beta<\kappa}cl_{\gamma}^{n}(G_{\beta})\neq X$ then $C_{\alpha}\setminus \bigcup _{\beta<\kappa}cl_{\gamma}^{n}(G_{\beta})\neq\emptyset$.
\end{enumerate}

Let $C=\bigcup_{\alpha<\kappa^{+}} C_{\alpha}$. We want to show that $X=C$. Assume there exists $y\in X\setminus C$. For every $B_{\beta}\in {\mathcal B}(y)$, $\beta\leq \kappa$, define ${\mathcal F}_{\beta}=\{V_{c}:\;c\in C\hbox{, and}\; cl_{\theta}^{n}(V_{c})\cap cl_{\theta}^{n}(B_{\beta})=\emptyset\}$. Since $X$ is $\theta^n$-Urysohn, we have $C\subseteq \bigcup_{\beta<\kappa}\bigcup {\mathcal F}_{\beta}$. By Lemma \ref{L1} we can find for every $\beta<\kappa$ a subcollection ${\mathcal G}_{\beta}\subseteq {\mathcal F}_{\beta}$, $|{\mathcal G}_{\beta}|\leq \kappa$ such that ${\mathcal F}_{\beta}\subseteq cl_{\gamma}^{n}(\bigcup_{G\in{\mathcal G}_{\beta}}cl_{\theta}^{n}(G))$. Note that $y\notin cl_{\gamma}^{n}(\bigcup_{G\in{\mathcal G}_{\beta}}cl_{\theta}^{n}(G))$ since $(\bigcup_{G\in {\mathcal G}_{\beta}}cl_{\theta}^{n}(G))\cap cl_{\theta}^{n}(B_{\beta})=\emptyset$. Find $\alpha<\kappa^{+}$ such that $\bigcup_{\beta<\kappa}{\mathcal G}_{\beta}\subseteq{\mathcal V}_{\alpha}$. Then $y\notin \bigcup_{\beta<\kappa}cl_{\gamma}^{n}(\bigcup_{G\in{\mathcal G}_{\beta}}cl_{\theta}^{n}(G))$ but $C_{\alpha}\subseteq C\subseteq \bigcup_{\beta<\kappa}\bigcup{\mathcal F}_{\beta}\subseteq \bigcup_{\beta<\kappa}cl_{\gamma}^{n}(\bigcup_{G\in{\mathcal G}_{\beta}}cl_{\theta}^{n}(G))$ and this is a contradiction.
\endproof

If $n=1$ we have the Schr\"oder inequality \cite{Sc}:
\begin{cor}
If $X$ is a Urysohn space, then $|X|\leq 2^{Uc(X)\chi(X)}$.
\end{cor}

\section{Cardinality bounds for homogeneous $n$-$\theta$-Urysohn spaces}\label{S3}

Many cardinality bounds for general spaces have corresponding ``companion'' cardinality bounds for homogeneous topological spaces. The latter utilizes homogeneity to give an improved bound. In Theorem~\ref{homogeneousbound} below we give the homogeneous companion bound to Theorem~\ref{cardbound}. Theorem~\ref{PHbound} generalizes this further to the power homogeneous setting. We recall the following definitions.

\begin{defin}
A space $X$ is \emph{homogeneous} if for all $x,y\in X$ there exists a homeomorphism $h:X\to X$ such that $h(x)=y$. $X$ is \emph{power homogeneous} if there exists a cardinal $\kappa$ such that $X^\kappa$ is homogeneous.
\end{defin}

The following theorem uses the Erd\"os-Rado Theorem and represents a variation of Proposition 2.1 in \cite{CR}.

\begin{theorem}\label{homogeneousbound}
If $X$ is homogeneous and $\theta^n$-Urysohn then $|X|\leq 2^{\theta^n-Uc(X)\pi_{\chi}(X)}$.
\end{theorem}

\proof
Let $\kappa=\theta^n$-$Uc(X)\pi_{\chi}(X)$. Fix a point $p\in X$ and a local $\pi$-base ${\mathcal B}$ at $p$ consisting of non-empty sets such that $|{\mathcal B}|\leq\kappa$. For all $x\in X$ let $h_x:X\rightarrow X$ be a homeomorphism such that $h_x(p)=x$. Let $\Delta=\{(x,x)\in[X]^2:x\in X\}$. Define $B:[X]^2\setminus\Delta\rightarrow{\mathcal B}$ as follows. For all $x\neq y$, there exist disjoint open sets $U(x,y)$ and $V(x,y)$ containing $x$ and $y$ respectively such that $cl_\theta^n(U(x,y))\cap cl_\theta^n(V(x,y))=\emptyset$. For each $x\neq y\in X$ the open set $h_x^{\leftarrow}[U]\cap h_y^{\leftarrow}[V]$ contains $p$. Thus there exists $B(x,y)\in{\mathcal B}$ such that $B(x,y)\subseteq h_x^{\leftarrow}[U]\cap h_y^{\leftarrow}[V]$. Observe that $h_x[cl_\theta^n(B(x,y))]\cap h_y[cl_\theta^n(B(x,y))]=\emptyset$ for all $(x,y)\in [X]^2\setminus\Delta$.

By way of contradiction suppose that $|X|>2^\kappa$. By the Erd\"os-Rado Theorem there exists $Y\in[X]^{\kappa^{+}}$ and $B\in{\mathcal B}$ such that $B=B(x,y)$ for all $x\neq y\in Y$. For $x\neq y\in Y$ we have $h_x[cl_\theta^n(B)]\cap h_y[cl_\theta^n(B)]=h_x[cl^n_\theta(B(x,y))]\cap h_y[cl^n_\theta(B(x,y))]=\emptyset$. This shows that ${\mathcal C}=\{h_x[cl_\theta^n(B)]:x\in Y\}$ is a $\theta^n$-Urysohn cellular family. But $|{\mathcal C}|=|Y|=\kappa^+>\theta^n\hbox{-}Uc(X)$, a contradiction. Therefore $|X|\leq 2^\kappa$.
\endproof

In the power homogeneous case, the proof is the same as the proof of Theorem 15 in \cite{BCCS}, except that closures are replaced with $\theta$-$n$-closures. In particular, this is done in the claim in that theorem.  

\begin{theorem}\label{PHbound}
If $X$ is power homogeneous and $\theta^n$-Urysohn then $|X|\leq 2^{\theta^n\hbox{-}Uc(X)\pi_{\chi}(X)}$.
\end{theorem}

\end{document}